\documentclass{article} 
\usepackage[a4paper]{anysize}
\usepackage{amsmath,amsfonts,amssymb,amsthm}
\usepackage{graphicx,pgf,tikz}
\usepackage{multirow}
\usepackage{float}
\usepackage{cite}

\usepackage[OT2,OT1]{fontenc}

\newtheorem{Theorem}{Theorem}
\newtheorem*{Remark}{Remark}

\title{On One Generalization of the Multipoint Nonlocal Contact
  Problem for Elliptic Equation in Rectangular Area}

\author{%
  Tinatin Davitashvili\footnote{Iv.Javakhishvili Tbilisi State
  University, FENS; tinatin.davitashvili@tsu.ge} \and
  Hamlet Meladze\footnote{Muskhelishvili Institute of Computational
  Mathematics, Georgian Technical University;
  h\_meladze@hotmail.com}\and
  Francisco Criado-Aldeanueva\footnote{Department of Applied Physics
  II, Polytechnic School, Malaga University;
  fcaldeanueva@ctima.uma.es} \and
  Jose Maria Sanchez\footnote{Department of Didactic of Mathematics,
  Education Faculty, Malaga University; jmss@uma.es} }

\begin{document}

\maketitle

\begin{abstract}
  A nonlocal contact problem for two-dimensional linear elliptic
  equations is stated and investigated. The method of separation of
  variables is used to find the solution of a stated problem in case
  of Poisson's equation. Then the more general problem with nonlocal
  multipoint contact conditions for elliptic equation with variable
  coefficients is considered and the iterative method to solve the
  problem numerically is constructed and investigated. The uniqueness
  and existence of the regular solution is proved. The iterative
  method allows to reduce the solution of a nonlocal contact problem
  to the solution of a sequence of classical boundary value problems.


\end{abstract}

\section{Introduction}\label{sec:1}

It can be stated that the study of non-local boundary and
initial-boundary problems and the development and analysis of
numerical methods for their solution is an important area of applied
mathematics. Nonlocal problems are naturally obtained in mathematical
models of real processes and phenomena in biology, physics,
engineering, ecology, etc. One can get acquainted with related
questions in works \cite{1,2,3} and the references cited there.

The publication of the articles \cite{4,5,6,7} laid the foundation of
further research in the area of non-local boundary problems and their
numerical solution. In 1969 was published the work of A.Bitsadze and
A.Samarskii \cite{7}, which was related to the mathematical modeling
of plasma processes. In this paper a new type of nonlocal problem for
elliptic equations was considered, hereinafter referred to as the
Bitsadze-Samarsky problem. But the intensive research into nonlocal
boundary value problems began in the 80s of the 20th century (see, for
instance, \cite{8,9,10,11,12,13,14,15,16,17,18,19,20,21} and
references herein).

The work \cite{8} is devoted to the formulation and investigation of a
non-local contact problem for a parabolic-type linear differential
equation with partial derivatives. In the first part of the work, the
linear parabolic equation with constant coefficients is considered. To
solve a non-local contact problem, the variable separation method
(Fourier method) is used. Analytic solutions are built for this
problem. Using the iterative method, the existence and uniqueness of
the classical solution to the problem is proved. The effectiveness of
the method is confirmed by numerical calculations. In paper \cite{9},
based on the variation approach, the definition of a classical
solution is generalized for the Bitsadze-Samarskii non-local boundary
value problem, posed in a rectangular area. In \cite{10} the
Bitsadze-Samarskii nonlocal boundary value problem for the
two-dimensional Poisson equation is considered for a rectangular
domain. The solution of this problem is defined as a solution of the
local Dirichlet boundary value problem, by constructing a special
method to find a function as the boundary value on the side of the
rectangle, where the nonlocal condition was given. In paper \cite{11}
the two-dimensional Poisson equation with nonlocal integral boundary
conditions in one of the directions, for a rectangular domain, is
considered. For this problem, a difference scheme of increased order
of approximation is constructed, its solvability is studied, and an
iterative method for solving the corresponding system of difference
equations is justified. In paper \cite{12} a two-step difference
scheme with second order of accuracy for an approximate solution of
the nonlocal boundary value problem for the elliptic differential
equation in an arbitrary Banach space with the strongly positive
operator is considered. In paper \cite{13} the optimal control problem
for Helmholtz equation with non-local boundary conditions and
quadratic functional is considered. The necessary and sufficient
optimal conditions in a maximum principle form have been obtained. In
\cite{14}, the Bitsadze-Samarskii boundary value problem is considered
for a linear differential equation of first order for the bounded
domain of the complex plane. The existence of a generalized equation
is proved and an a priori estimate is obtained. Then the corresponding
theorem on existence and uniqueness of a generalized solution is
proved. A boundary-value problem with a nonlocal integral condition is
considered in \cite{15} for a two-dimensional elliptic equation with
mixed derivatives and constant coefficients. The existence and
uniqueness of a weak solution is proved in a weighted Sobolev space. A
difference scheme is constructed and its convergence is proved. In
paper \cite{16} the one class of nonlocal in time problems for
first-order evolution equations is considered. The solvability of the
stated problem is investigated. All of them are, basically, related to
the problems with nonlocal conditions considered only at the border of
the area of definition of the differential operator.

In the present paper, the multipoint nonlocal contact problem for
linear elliptic equations is stated and investigated in
two-dimensional domains. The method of separation of variables is used
to find the solution of a stated problem in case of Poisson's
equation. Then the more general problem with nonlocal contact
conditions for elliptic equation with variable coefficients is
considered and the iterative method to numerically solve the problem
is constructed and investigated. The iterative method allows to reduce
the solution of a nonlocal contact problem to the solution of a
sequence of classical boundary problems. The numerical experiment is
conducted. The results fully agree with the theoretical conclusions
and show the efficiency of the proposed iterative procedure.

\section{Method of Separation of Variables for Poisson Equation}\label{sec:2}

\subsection{Formulation of the problem}\label{sec:2.1}
\begin{figure}[htbp]
  \centering
  \includegraphics{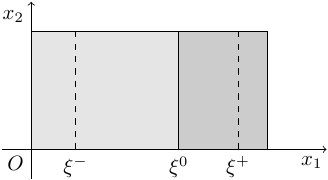}
  \caption{Domain for Poisson equation}
  \label{fig:1}
\end{figure}

Let us consider rectangular domain in two-dimensional space $R^2$ with
boundary $\Gamma$:
\begin{equation*}
  \{(x_1,x_2 )\ |\  0\le x_1\le 1,\ 0\le x_2\le 1\}.
\end{equation*}
Suppose, $0 < \xi^0 < 1$, and define the segment
$\{(x_1,x_2)\ |\ x_1=\xi^0,\ 0 \le x_2 \le 1\}$ (see Figure
\ref{fig:1}).

We consider the following nonlocal contact problem: Find the
continuous function
\begin{equation}\label{eq:1}
  u(x_1,x_2 )=
  \begin{cases}
    u^- (x_1,x_2 ), & 0\le x_1 < \xi^0,\ 0\le x_2\le 1,\\
    u(\xi^0,x_2 ),  & x_1=\xi^0,\ 0\le x_2\le 1,\\ 
    u^+ (x_1,x_2 ), & \xi^0 < x_1\le 1,\ 0\le x_2\le 1,
  \end{cases}  
\end{equation}
which satisfies the equations:
\begin{equation}\label{eq:2}
  \begin{aligned}
    \Delta u^- (x_1,x_2 )&=f^- (x_1,x_2 ),&& 0< x_1< \xi^0,\ 0< x_2< 1, \\
    \Delta u^+ (x_1,x_2 )&=f^+ (x_1,x_2 ),&& \xi^0<x_1<1,\ 0<x_2<1,
  \end{aligned}
\end{equation}
the boundary conditions
\begin{equation}
  \begin{gathered}
    \left.
    \begin{aligned}
      u^- (x_1,0)&=0, &&0\le x_1\le \xi^0,\\
      u^+ (x_1,0)&=0, &&\xi^0\le x_1\le 1,
    \end{aligned}
    \right\}
    \\
    \left.\begin{aligned}
      u^- (x_1,1)&=0, && 0\le x_1\le \xi^0,\\
      u^+ (x_1,1)&=0, && \xi^0\le x_1\le 1,    
    \end{aligned}\right\}
    \\
    \left.\begin{aligned}
      u^- (0,x_2)&=0, && 0\le x_2\le 1,\\
      u^+ (1,x_2)&=0, && 0\le x_2\le 1, 
    \end{aligned}\right\}
  \end{gathered}\label{eq:5}
\end{equation}
and the nonlocal contact condition
\begin{multline}\label{eq:6}
  u^-(\xi^0,x_2)=u^+ (\xi^0,x_2 )=u(\xi^0,x_2)=u_0(x_2)\\ 
  =\gamma^- u^-(\xi^-,x_2)+\gamma^+ u^+(\xi^+,x_2)+\phi_0(\xi^0,x_2 ), \qquad 0 < x_2 < 1,
\end{multline}
where $f^- (x_1,x_2 ),\, f^+ (x_1,x_2 )$ are known, sufficiently
smooth functions, and
\begin{equation}\label{eq:7}
  0<\xi^-<\xi^0,\,
  \xi^0<\xi^+<1,\,
  \gamma^->0,\,
  \gamma^+>0,\,
  \gamma^-+\gamma^+\le 1.
\end{equation}

Note that in previously published articles the nonlocal conditions
were mainly formulated under the restriction of the following type:
$\gamma^- + \gamma^+ < 1$. In this article the results are achieved
considering the following conditions $\gamma^-+\gamma^+ \le 1$.

\subsection{Separation of variables}\label{sec:2.2}

Using the method of separation of variables, we can build the solution
of nonlocal contact problem (\ref{eq:1})-(\ref{eq:7}). Note that this
technique can be extended for more general case.

We will find the solution of nonlocal contact problem
(\ref{eq:1})-(\ref{eq:7}) in a following form:
\begin{align}
  u^-(x_1,x_2 )&=\sum_{k=1}^{\infty}a_k^- (x_1 )\sin k\pi x_2, && 0\le x_1\le \xi^0,\ 0\le x_2\le 1,\label{eq:8}\\
  u^+(x_1,x_2 )&=\sum_{k=1}^{\infty}a_k^+ (x_1 )\sin k\pi x_2, && \xi^0\le x_1\le 1,\ 0\le x_2\le 1.\label{eq:9}      
\end{align}

It is evident that the functions (\ref{eq:8}) and (\ref{eq:9}) satisfy
the boundary conditions:
\begin{align*}
  u^- (x_1,0)&=u^- (x_1,1)=0,&& 0\le x_1\le \xi^0,\\
  u^+ (x_1,0)&=u^+ (x_1,1)=0,&& \xi^0\le x_1\le 1.
\end{align*}

We must choose the coefficients $a_k^+(x_1)$ and
$a_k^-(x_1),\ k=1,2,3,\ldots,$ so that the functions (\ref{eq:8}),
(\ref{eq:9}) satisfy the equations (\ref{eq:2}), the rest boundary
conditions
\begin{equation*}
  u^-(0,x_2)=0,\ u^+(1,x_2)=0,\qquad 0\le x_2\le 1,
\end{equation*}
and the nonlocal contact condition as well. Thus, $a_k^+ (x_1 )$ and
$a_k^- (x_1 ),\ k=1,2,3,…,$ should be solutions of the following
problems:
\begin{align}
  \frac{d^2 a_k^-(x_1)}{dx_1^2} -\lambda_k^2 a_k^- (x_1 ) & =f_k^- (x_1),
  &&  0<x_1<\xi^0,\  a_k^-(0)=0,\ a_k^-(\xi^0)=\Phi_k, \label{eq:10}\\
  \frac{d^2 a_k^+(x_1)}{dx_1^2} -\lambda_k^2 a_k^+ (x_1 ) &=f_k^+ (x_1),
  && \xi^0<x_1<1,\  a_k^+(\xi^0)=\Phi_k,\ a_k^+(1)=0,\label{eq:11}           
\end{align} 
where $\lambda_k=\pi^2 k^2,\ f_k^- (x_1 ),\ f_k^+ (x_1 )$ are
the coefficients of Fourier Series expansion of the functions $f^-
(x_1,x_2 )$ and $f^+ (x_1,x_2 )$:
\begin{equation*}
  f^-(x_1,x_2 )=\sum_{k=1}^{\infty}f_k^-(x_1 ) \sin \lambda_k x_2,\qquad
  f^+(x_1,x_2 )=\sum_{k=1}^{\infty}f_k^+ (x_1 ) \sin \lambda_k x_2,
\end{equation*}
and $\Phi_k$, $k=1,2,\ldots,$ are so far unknown constants
\begin{equation}\label{eq:12} 
u(\xi^0,x_2 )=u_0(x_2)=\sum_{k=1}^{\infty}\Phi_k \sin\lambda_k x_2.       
\end{equation}

\subsection{Existence and uniqueness of the solution}\label{sec:2.3}

Using the method of constructing general solutions of differential
equations with homogeneous boundary conditions and method of variation
of constants or method of undetermined coefficients of Lagrange (see,
for example \cite{22}), it is possible to build the general solutions
of equations (\ref{eq:10}) and (\ref{eq:11}).

At first, we will consider the problem (\ref{eq:10}). The general
solution of the problem (\ref{eq:10}) can be written in the following
form:
\begin{multline*}
  a_k^-(x_1)
  = \tilde{c}_1^-(x_1) e^{\sqrt{\lambda_k} x_1}
  + \tilde{c}_2^-(x_1) e^{-\sqrt{\lambda_k} x_1}\\
  + \frac{1}{\sqrt{\lambda_k}} \int_0^{x_1} \sinh\left(\sqrt{\lambda_k} (x_1-s)\right) f_k^-(s)\, ds,
  \qquad k=1,2,\ldots,
\end{multline*}
where we can define $\tilde{c}_1^-(x_1)$, $\tilde{c}_2^-(x_1)$,
considering the boundary conditions $a_k^-(0) = 0$, $a_k^-(\xi^0)=
\Phi_k$.

Finally the general solution of the problem (\ref{eq:10}) will get the
following form:
\begin{multline*}
  a_k^-(x_1) = \frac{(\Phi_k-I^-) \sinh(\sqrt{\lambda_k} x_1)}{\sinh(\sqrt{\lambda_k} \xi^0)}
  + \frac{1}{\sqrt{\lambda_k}} \int_0^{x_1}\sinh\left(\sqrt{\lambda_k}(x_1-s)\right) f_k^-(s)\,ds,\\
  k=1,2,\ldots
\end{multline*}
where
\begin{equation}
  I^- = \frac{1}{\sqrt{\lambda_k}} \int_0^{\xi^0} \sinh\left(\sqrt{\lambda_k} (\xi^0-s)\right)\, f_k^-(s)\,ds.
  \label{eq:15}
\end{equation}

Consequently, the formal solution of the problem
(\ref{eq:2})-(\ref{eq:5}) on the left subarea is the following
function
\begin{multline*}
  u^-(x_1,x_2)
  = \sum_{k=1}^\infty\left[
    \frac{(\Phi_k - I^-) \sinh(\sqrt{\lambda_k} x_1)}{\sinh(\sqrt{\lambda_k} \xi^0)}\right.\\
    \left.
    + \frac 1{\sqrt{\lambda_k}} \int_0^{x_1} \sinh\left(\sqrt{\lambda_k} (x_1-s)\right) f_k^-(s)\,ds
    \right]
  \sin \lambda_k x_2,
\end{multline*}
where $I^-$ is defined using (\ref{eq:15}).

Analogously, the general solution of the problem (\ref{eq:11}) can be
written in the following form:
\begin{multline*}
  a_k^+(x_1)
  = \tilde c_1^+(x_1) e^{\sqrt{\lambda_k} x_1}
  + \tilde c_2^+(x_1) e^{-\sqrt{\lambda_k} x_1}\\
  + \frac{1}{\sqrt{\lambda_k}} \int_{\xi^0}^{x_1} \sinh\left(\sqrt{\lambda_k} (x_1-s)\right) f_k^+(s)\,ds,
  \qquad k=1,2,\ldots
\end{multline*}
Considering the boundary conditions $a_k^+(\xi^0)=\Phi_k$,
$a_k^+(1)=0$, we can define $\tilde c_1^+$, $\tilde c_2^+$ uniquely.
  
Finally we can describe the formal solution of the problem
(\ref{eq:2})-(\ref{eq:5}) on the right subarea in a following way:

\begin{multline*}
  u^+(x_1,x_2 )=\sum_{k=1}^{\infty}\left[
    \frac{(\Phi_k\sinh(\sqrt{\lambda_k}(x_1-1))-I^+\sinh(\sqrt{\lambda_k}(x_1-\xi^0))}{\sinh(\sqrt{\lambda_k}(\xi^0-1))}\right.\\
    \left.+\frac{1}{\sqrt{\lambda_k }}\int_{\xi^0}^{x_1}\sinh(\sqrt{\lambda_k}(x_1-s)) f_k^+(s)ds\right]\sin\lambda_k x_2,
\end{multline*}
where 
\begin{equation}\label{eq:22}
  I^+=\frac{1}{\sqrt{\lambda_k}}\int_{\xi^0}^{1}\sinh(\sqrt{\lambda_k}(1-s)) f_k^+(s)ds.
\end{equation}

We can define the coefficients $\Phi_k$, $k=1,2,\ldots$, using the equality
(\ref{eq:12}) and nonlocal contact condition (\ref{eq:6}):
\begin{multline*}
  \Phi_k=\gamma^-\frac{\Phi_k \sinh(\sqrt{\lambda_k}\xi^-)}{\sinh(\sqrt{\lambda_k } \xi^0 )} +\gamma^+\frac{\Phi_k \sinh(\sqrt{\lambda_k}(\xi^+-1))}{\sinh(\sqrt{\lambda_k}(\xi^0-1))}\\
  -(F_k^-+F_k^+ )+
  \frac{\gamma^-}{\sqrt{\lambda_k }}\int_{0}^{\xi^-}\sinh(\sqrt{\lambda_k}(\xi^--s)) f_k^-(s)ds\\
  +\frac{\gamma^+}{\sqrt{\lambda_k}}\int_{\xi^0}^{\xi^+}\sinh(\sqrt{\lambda_k}(\xi^+-s)) f_k^+(s)ds+\phi_{0k}, 
\end{multline*}
where 
\begin{align*}
  \phi_0(\xi^0,x_2) &=\sum_{k=1}^{\infty}\phi_{0k} \sin \lambda_k x_2,\nonumber\\
  F_k^-&=\gamma^- \frac{\sinh(\sqrt{\lambda_k} \xi^-)}{\sqrt{\lambda_k} \sinh(\sqrt{\lambda_k}\xi^0)}\int_{0}^{\xi^0}\sinh(\sqrt{\lambda_k}(\xi^0-s)) f_k^-(s)ds,
  \\
  F_k^+ &=\gamma^+ \frac{\sinh(\sqrt{\lambda_k} (\xi^+ -\xi^0))}{\sqrt{\lambda_k} \sinh(\sqrt{\lambda_k}(\xi^0-1))}\int_{\xi^0}^{1}\sinh(\sqrt{\lambda_k}(1-s)) f_k^+(s)ds.
\end{align*}

Then we will get   
\begin{equation}\label{eq:26}
  \Phi_k\left\{1-\left[
    \gamma^- \frac{\sinh(\sqrt{\lambda_k} \xi^-)}{\sinh(\sqrt{\lambda_k}\xi^0)}
    +\gamma^+ \frac{\sinh(\sqrt{\lambda_k} (\xi^+ -1))}{\sinh(\sqrt{\lambda_k}(\xi^0-1))}\right] \right\}=-F_k+\phi_{0k},
\end{equation}
where
\begin{multline}\label{eq:27}
  F_k=F_k^- +F_k^+
  -\left[
    \frac{\gamma^-}{\sqrt{\lambda_k }}\int_{0}^{\xi^-}\sinh(\sqrt{\lambda_k}(\xi^--s)) f_k^-(s)ds\right.\\
    \left.+\frac{\gamma^+}{\sqrt{\lambda_k}}\int_{\xi^0}^{\xi^+}\sinh(\sqrt{\lambda_k}(\xi^+-s)) f_k^+(s)ds
  \right]. 
\end{multline}

As
\begin{equation*}
  \frac{\sinh(\sqrt{\lambda_k} \xi^-)}{\sinh(\sqrt{\lambda_k}\xi^0)}<1\  \mbox{  and  }\  \frac{\sinh(\sqrt{\lambda_k} (\xi^+ -1))}{\sinh(\sqrt{\lambda_k}(\xi^0-1))}<1,
\end{equation*}
then we will have
\begin{equation*}
  1-\left[
    \gamma^- \frac{\sinh(\sqrt{\lambda_k} \xi^-)}{\sinh(\sqrt{\lambda_k}\xi^0)}
    +\gamma^+ \frac{\sinh(\sqrt{\lambda_k} (\xi^+ -1))}{\sinh(\sqrt{\lambda_k}(\xi^0-1))}
    \right]>1-(\gamma^- + \gamma^+)\ge 0.
\end{equation*}

Consequently, from the equality (\ref{eq:26}) we get
\begin{multline}\label{eq:29}
  \Phi_k=\left\{1-
  \left[\gamma^- \frac{\sinh(\sqrt{\lambda_k} \xi^-)}{\sinh(\sqrt{\lambda_k}\xi^0)}
    +\gamma^+ \frac{\sinh(\sqrt{\lambda_k} (\xi^+ -1))}{\sinh(\sqrt{\lambda_k}(\xi^0-1))}
    \right]
  \right\}^{-1}(-F_k+\phi_{0k}), \\
  k=1,2,\ldots,
\end{multline}
where $F_k$ is defined from (\ref{eq:27}).  Finally, the formal
solution of the problem (\ref{eq:2})-(\ref{eq:6}) is the following
function:
\begin{equation*}
  u(x_1,x_2 )=
  \begin{cases}
    u^- (x_1,x_2 ), & 0\le x_1 < \xi^0,\ 0\le x_2\le 1,\\
    u(\xi^0,x_2 ), & x_1=\xi^0,\ 0\le x_2\le 1,\\
    u^+ (x_1,x_2 ), & \xi^0 < x_1\le 1,\ 0\le x_2\le 1,
  \end{cases}  
\end{equation*}
where 
\begin{multline*}
  u^- (x_1,x_2 )=\sum_{k=1}^{\infty}\left[
    \frac{(\Phi_k-I^-)\sinh(\sqrt{\lambda_k} x_1)}{\sinh(\sqrt{\lambda_k}\xi^0)}\right.\\
    \left.+\frac{1}{\sqrt{\lambda_k }}\int_{0}^{x_1}\sinh(\sqrt{\lambda_k}(x_1-s)) f_k^-(s)ds
    \right]\sin\lambda_k x_2,
\end{multline*}
\begin{equation*}
  u(\xi^0,x_2 )=\sum_{k=1}^{\infty}\Phi_k \sin\lambda_k x_2,
\end{equation*}
\begin{multline*}
  u^+ (x_1,x_2 )=\sum_{k=1}^{\infty}\left[
    \frac{(\Phi_k\sinh(\sqrt{\lambda_k}(x_1-1))-I^+\sinh(\sqrt{\lambda_k}(x_1-\xi^0))}{\sinh(\sqrt{\lambda_k}(\xi^0-1))}\right.\\
    \left.+\frac{1}{\sqrt{\lambda_k }}\int_{\xi^0}^{x_1}\sinh(\sqrt{\lambda_k}(x_1-s)) f_k^+(s)ds\right]\sin\lambda_k x_2, 
\end{multline*}
$\Phi_k$ is defined from (\ref{eq:29}), and $I^-,\ I^+$ - from
(\ref{eq:15}) and (\ref{eq:22}), respectively.

Thus, the following theorem is true.

\begin{Theorem}
  If $f^- (x_1,x_2 ),\ f^+ (x_1,x_2 )$ and $\phi_0(\xi^0,x_2 )$ are
  sufficiently smooth functions, then the problem
  (\ref{eq:2})-(\ref{eq:6}) has a unique regular solution.
\end{Theorem}

Note, that the applied technique can be successfully used for more
general problems, but in this case the use of spectral theory of
linear operators will be needed.


\section{Nonlocal Contact Problem for Equation with Variable Coefficients}\label{sec:3}

\subsection{Designations}\label{sec:3.1}

\begin{figure}[htbp]
  \centering
  \includegraphics{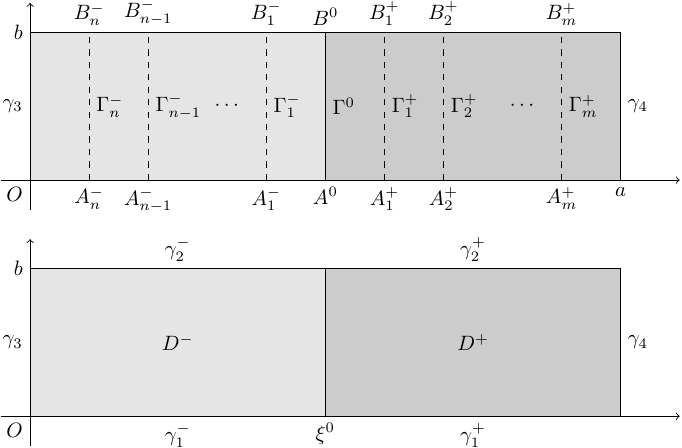}
  \caption{Domain for contact problem for equation with variable
    coefficients, $\bar D = D\cup\gamma$}
  \label{fig:2}
\end{figure}

Now let us consider the problem with nonlocal contact conditions for
elliptic equation with variable coefficients.

Suppose $D$ is a rectangular domain in two-dimensional space $R^2$
(see Figure \ref{fig:2}): $D=\{(x_1,x_2 )\ |\ 0<x_1<a,\ 0<x_2<b\}$
with piecewise boundary $\gamma=\bigcup_{i=1}^4\gamma_i$, where
\begin{align*}
  \gamma_1 &=\{(x_1,x_2 )\ |\ 0\le x_1\le a,\ x_2=0\}, &
  \gamma_2&=\{(x_1,x_2 )\ |\ 0\le x_1\le a,\ x_2=b\},\\
  \gamma_3&=\{(x_1,x_2 )\ |\ x_1=0,\ 0\le x_2\le b\}, &
  \gamma_4&=\{(x_1,x_2 )\ |\ x_1=a,\  0\le x_2\le b\}.
\end{align*}

Suppose, $0 < \xi^0 < a$, and define:
\begin{align*}
  \Gamma^0 &= \{(x_1 , x_2 )\ |\ x_1 = \xi^0 , 0 \le x_2 \le b\}\\
  \gamma_1^- &= \{(x_1 , x_2)\ |\ 0 \le x_1 \le \xi^0 , x_2 = 0\}, &
  \gamma_2^- &= \{(x_1 , x_2)\ |\ 0 \le x_1 \le \xi^0 , x_2 = b\},\\
  \gamma_1^+ &= \{(x_1 , x_2 )\ |\ \xi^0 \le x_1 \le a, x_2 = 0\}, &
  \gamma_2^+ &= \{(x_1 , x_2 )\ |\ \xi^0 \le x_1 \le a, x_2 = b\}.
\end{align*}
It is obvious, that $\gamma_1^- \cup \gamma_1^+ = \gamma_1$,
$\gamma_2^- \cup \gamma_2^+ = \gamma_2$, and $\Gamma^0$ divides the
domain $D$ into two parts (subdomains) $D^-$ and $D^+$ , where
\begin{equation*}
  D^- = \{(x_1 , x_2 )\ |\ 0 < x_1 < \xi^0 , 0 < x_2 < b\}, \qquad
  D^+ = \{(x_1 , x_2 )\ |\ \xi^0 < x_1 < a, 0 < x_2 < b\}
\end{equation*}

\begin{multline*}
  \Gamma_i^-=\{(x_1,x_2)\ |\ x_1=\xi_i^-,\ 0<\xi_i^-<\xi^0,\ 0\le x_2\le b\},\\
  i=¯\overline{(1,n)},\ 0<\xi_n^-<\cdots<\xi_1^-<\xi^0,
\end{multline*}
\begin{multline*}
  \Gamma_j^+=\{(x_1,x_2 )\ |\ x_1=\xi_j^+,\ \xi^0<\xi_j^+<a,\ 0\le x_2\le b\},\\
  j=\overline{(1,m)},
  \ \xi^0<\xi_1^+<\xi_2^+<\cdots<\xi_m^+<b.
\end{multline*}
$\Gamma^0, \Gamma_i^-,\ i=\overline{(1,n)},$ and $\Gamma
_j^+,\ j=\overline{(1,m)},$ intersect $\gamma_1$ and $\gamma_2$,
respectively, in the points:
\begin{equation*}
  A^0(\xi^0,0),\
  B^0(\xi^0,b),\
  A_i^-(\xi_i^-,0),\
  B_i^-(\xi_i^-,b),\
  A_i^+(\xi_i^+,0)  \text{ and }
  B_i^+(\xi_i^+,b) .
\end{equation*}

%

\subsection{Statement of the Problem}\label{sec:3.2}

Let us consider the following problem: Find in domain $\bar D = D \cup
\gamma$ (where $\gamma$ is defined in \ref{sec:3.1}) a continuous
function $u(x_1, x_2)$:
%
\begin{equation}\label{eq:31}
  \begin{gathered}
    u(x_1,x_2 )=
    \begin{cases}
      u^- (x_1,x_2 ),&\mbox{ if } (x_1,x_2)\in D^-,\\
      u^0(x_1,x_2 ),&\mbox{ if } (x_1,x_2)\in \Gamma^0,\\
      u^+ (x_1,x_2 ),&\mbox{ if } (x_1,x_2)\in D^+,
    \end{cases}\\
    u^-(x_1,x_2 )\in C^2(D^-),\
    u^+(x_1,x_2 )\in C^2(D^+),\
    u^0(x_1,x_2)\in C(\Gamma ^0),
  \end{gathered}
\end{equation}
which satisfies the equations 
\begin{multline}\label{eq:32} 
  L^- u^-(x_1,x_2 )\equiv\\
  \sum_{\alpha,\beta=1}^2 \frac{\partial}{\partial x_\alpha}(K_{\alpha \beta}^- (x_1,x_2 )\frac{\partial u^-}{\partial x_\beta})-k^- (x_1,x_2 ) u^-
  =-f^-(x_1,x_2 ),\  (x_1,x_2 )\in D^-, 
\end{multline}
\begin{multline}\label{eq:33} 
  L^+ u^+(x_1,x_2 )\equiv\\
  \sum_{\alpha,\beta=1}^2 \frac{\partial}{\partial x_\alpha} (K_{\alpha \beta}^+ (x_1,x_2)\frac{\partial u^+}{\partial x_\beta })-k^+(x_1,x_2 )u^+=-f^+ (x_1,x_2 ),\ (x_1,x_2 )\in  D^+, 
\end{multline}

The function $u(x_1,x_2 )$ also satisfies the boundary conditions
\begin{align*}
  u^-(x_1,x_2 )&=\phi^- (x_1,x_2 ), &&  (x_1,x_2 )\in \gamma_1^-\cup \gamma_2^-\cup \gamma_3,
  \\
  u^+(x_1,x_2)&=\phi^+(x_1,x_2 ), && (x_1,x_2)\in \gamma_1^+\cup \gamma_2^+\cup \gamma_4, 
\end{align*}
the nonlocal contact conditions
\begin{equation*}
  u^-(\Gamma^0)=u^+(\Gamma^0)=u(\Gamma^0)=\sum_{i=1}^n \beta_i^- u^-(\Gamma_i^-)+\sum_{j=1}^m \beta_j^+  u^+(\Gamma_j^+ )+\phi^0 (\Gamma^0),     
\end{equation*}
and the coordination conditions 
\begin{equation}\label{eq:37}
  \begin{aligned}
    u(A^0)&=\sum_{i=1}^n \beta_i^- u^-(A _i^-)+\sum_{j=1}^m \beta_j^+  u^+(A_j^+ )+\phi^0 (A^0),\\
    u(B^0)&=\sum_{i=1}^n \beta_i^- u^-(B_i^-)+\sum_{j=1}^m \beta_j^+  u^+(B_j^+ )+\phi^0 (B^0),   
\end{aligned}       
\end{equation}
where $\beta_i^-=\text{const} > 0$, $\beta_j^+=\text{const} > 0$,
$0<\sum_{i=1}^n \beta_i^- +\sum_{j=1}^m \beta_j^+ \le 1.$ $K_{\alpha
  \beta}^\pm (\cdot),$ $k^\mp (\cdot)\ f^\pm (\cdot)$, $\phi^\pm
(\cdot)$ and $\phi^0 (\cdot)$ are known functions, which satisfy all
the conditions, that provide existence of the unique solutions of
Dirichlet problem in $D^-$ and $D^+$ \cite{23}.

Suppose that the equations (\ref{eq:32}) and (\ref{eq:33}) are
uniformly elliptic. Then their coefficients satisfy the following
conditions \cite{23}:
%
%
%
\begin{equation*}
\begin{aligned}
  4K_{11}^- (x_1,x_2 ) \cdot K_{22}^- (x_1,x_2 )&>
  (K_{12}^- (x_1,x_2 )+K_{21}^- (x_1,x_2 ))^2,\\
  4K_{11}^+ (x_1,x_2 ) \cdot K_{22}^+ (x_1,x_2 )&>
  (K_{12}^+ (x_1,x_2 )+K_{21}^+ (x_1,x_2 ))^2.         
\end{aligned}       
\end{equation*}
$K^-_{\alpha\beta}(x_1,x_2)$ and $u^-(x_1,x_2)$ can be considered as
coefficient of heat conductivity and temperature of the first body
($D^-$), $K^+_{\alpha\beta}(x_1,x_2)$ and $u^+(x_1,x_2)$ - of the
second body ($D^+$). Thus, the stated problem can be considered as a
mathematical model of stationary distribution of heat in two
contacting isotropic bodies.

We will call the problem (\ref{eq:31})-(\ref{eq:37}) nonlocal contact
one since it is generalization of a classical contact problem.

\subsection{Uniqueness of a Solution of Problem (\ref{eq:31})-(\ref{eq:37})}\label{sec:3.3}

The following theorem is true.

\begin{Theorem}\label{thm:2}
  If the regular solution of problem (\ref{eq:31})-(\ref{eq:37})
  exists and condition $0<\sum_{i=1}^n \beta_i^- +\sum_{j=1}^m
  \beta_j^+ \le 1$ is fulfilled, then then the solution is unique.\\
\end{Theorem}

\begin{proof}
  Suppose that problem (\ref{eq:31})-(\ref{eq:37}) has two solutions:
  $v(x_1,x_2 )$ and $w(x_1,x_2 )$. Then for the function
  $z(x_1,x_2)=v(x_1,x_2 )-w(x_1,x_2 )$ we will have the following
  problem
  \begin{align*}
    L^- z^- (x_1,x_2 )&=0, \text{ if }  (x_1,x_2 )\in D^-, 
    \\
    L^+ z^+ (x_1,x_2 )&=0, \text{ if } (x_1,x_2 )\in D^+, 
  \end{align*}
  \begin{equation}\label{eq:41}
    \begin{aligned}
      z^- (x_1,x_2 )&=0, &&\text{ if } (x_1,x_2 )\in \gamma_1^-\cup \gamma_2^-\cup \gamma_3,\\
      z^+ (x_1,x_2 )&=0, &&\text{ if }  (x_1,x_2 )\in \gamma_1^+\cup \gamma_2^+\cup \gamma_4,
    \end{aligned}
  \end{equation}
  \begin{equation}\label{eq:42}
    z(\Gamma^0)=z^- (\Gamma^0)=z^+ (\Gamma^0)=\sum_{i=1}^n \beta_i^-  z^-(\Gamma_i^-)+\sum_{j=1}^m \beta_j^+  z^+(\Gamma_j^+).     
  \end{equation}

  From the equality (\ref{eq:42}) it follows that
  
\begin{multline*}
    \max\vert z(\Gamma^0)\vert\le \max \sum_{i=1}^n \beta_i^- \vert z^-(\Gamma_i^-)\vert +\max \sum_{j=1}^m \beta_j^+\vert z^+ (\Gamma_j^+)\vert \\
    \le \max_{0\le i \le n} \vert z^-(\Gamma_i^-)\vert \sum _{i=1}^n \beta_i^- +\max_{0\le j \le m} |z^+ (\Gamma_j^+)\vert \sum_{j=1}^m \beta_j^+. 
  \end{multline*}

  Taking into account the condition $0<\sum_{i=1}^n \beta_i^-
  +\sum_{j=1}^m \beta_j^+ \le 1$, we obtain
  
\begin{equation*}
    \max\vert z (\Gamma_0)\vert \le \max_{0\le i \le n} \vert z^-(\Gamma_i^-)\vert \ \mbox{  or }\  \max\vert z (\Gamma_0)\vert \le\max_{0\le j \le m} |z^+ (\Gamma_j^+)\vert.
  \end{equation*}

  This means that the function z does not attain a maximum on
  $\Gamma_0$, but attains its maximum on $D^-$ or $D^+$, that
  contradicts the maximum principle. So, $z\equiv \text{const}$ and
  taking into account condition (\ref{eq:41}), we obtain $z\equiv0$,
  i.e. the solution of the problem (\ref{eq:31}) - (\ref{eq:37}) is
  unique.
\end{proof}

\subsection{Iterative Method for Problem (\ref{eq:31})-(\ref{eq:37})}\label{sec:3.4}

Let us consider the following iteration process for the problem
(\ref{eq:31})-(\ref{eq:37}):
\begin{align}
  L^-[u^-(x_1,x_2 )]^{(k)}&= -f^-(x_1,x_2 ),&& \mbox{if } (x_1,x_2 )\in D^-,\label{eq:43}\\
  L^+[u^+(x_1,x_2 )]^{(k)}&= -f^+ (x_1,x_2 ),&& \mbox{if }(x_1,x_2 )\in  D^+, \label{eq:44} \\
  [u^-(x_1,x_2)]^{(k)}&= \phi^- (x_1,x_2 ), &&  \mbox{if }(x_1,x_2 )\in \gamma_1^-\cup \gamma_2^-\cup \gamma_3,\label{eq:45} \\
  [u^+(x_1,x_2)]^{(k)}&=\phi^+(x_1,x_2), && \mbox{if }(x_1,x_2)\in \gamma_1^+\cup \gamma_2^+\cup \gamma_4,\label{eq:46} 
\end{align}
\begin{multline}\label{eq:47} 
  u^{(k)}(\Gamma^0)=[u^-(\Gamma^0)]^{(k)}=[u^+(\Gamma^0)]^{(k)}\\
  =\sum_{i=1}^n \beta_i^- [u^-(\Gamma_i^-)]^{(k-1)}+\sum_{j=1}^m \beta_j^+[u^+(\Gamma_j^+)]^{(k-1)}+\phi^0 (\Gamma^0),     
\end{multline}
where $k=1,2,\ldots$. Initially we can take e.g.
$[u^-(\Gamma_i^-)]^{(0)}=0,$
$[u^+(\Gamma_j^+)]^{(0)}=0,$
$i=\overline{(1,n)},\ j=\overline{(1,m)}.$

Given the initial approximations in nonlocal contact condition
(\ref{eq:47}) of the iterative process (\ref{eq:41})-(\ref{eq:47}),
$\left[u^-(\Gamma^-_i)\right]^{(k-1)}$ and
$\left[u^+(\Gamma^+_j)\right]^{(k-1)}$, $i = \overline{(1, n)}$, $j =
\overline{(1, m)}$, we can calculate the values of $u$ on $\Gamma^0$
and, thus, get two classical boundary problems. After solving these
problems, we can define the consequent values of $u$ on $\Gamma^0$
from (\ref{eq:47}) for the next iteration, etc.

\begin{Theorem}\label{thm:3}
  If the solution of problem (\ref{eq:31})-(\ref{eq:37}) exists, then
  the iterative process (\ref{eq:43})-(\ref{eq:47}) converges to this
  solution at the rate of an infinitely decreasing geometric
  progression.
\end{Theorem}

\begin{proof}
  Denote by $z^{(k)}(x_1,x_2) = u^{(k)}(x_1,x_2) - u(x_1,x_2)$, where
  $u$ is a solution of the problem (\ref{eq:31})-(\ref{eq:37}) and
  $u^{(k)}$ -- of the problem (\ref{eq:43})-(\ref{eq:47}).
    
  Then we obtain the following problems:
  %
  \begin{align*}
    L^-[z^-(x_1,x_2 )]^{(k)}&= 0,&& \text{if } (x_1,x_2 )\in D^-, 
    \\
    L^+[z^+(x_1,x_2 )]^{(k)}&= 0,&& \text{if }(x_1,x_2 )\in  D^+, 
    \\
    [z^-(x_1,x_2)]^{(k)} &= 0,  &&\text{if }(x_1,x_2 )\in \gamma_1^-\cup \gamma_2^-\cup \gamma_3, 
    \\
    [z^+(x_1,x_2)]^{(k)}&=0, && \text{if }(x_1,x_2)\in \gamma_1^+\cup \gamma_2^+\cup \gamma_4, 
  \end{align*}
  \begin{multline}\label{eq:52} 
  [z(\Gamma^0)]^{(k)}=[z^-(\Gamma^0)]^{(k)}=[z^+(\Gamma^0)]^{(k)} \\
  =\sum_{i=1}^n \beta_i^- [z^-(\Gamma_i^-)]^{(k-1)}+\sum_{j=1}^m \beta_j^+[z^+(\Gamma_j^+)]^{(k-1)},     
\end{multline}
where $k=1,2,\ldots,$ and $[z^-(\Gamma_i^-)]^{(0)}=0,$ $[z^+
  (\Gamma_j^+)]^{(0)}=0,$
$i=\overline{(1,n)},\ j=\overline{(1,m)}.$\\ From the equality
(\ref{eq:52}) we have
\begin{equation*}
  \max_{\Gamma^0}\left|[z(\Gamma^0)]^{(k)} \right| \le
  \max_{1 \le i \le n}\left| [z^-(\Gamma_i^-)]^{(k-1)}\right|\sum_{i=1}^n \beta_i^-
  +\max_{1 \le j \le m} \left| [z^+(\Gamma_j^+)]^{(k-1)}\right|\sum_{j=1}^m \beta_j^+.
\end{equation*}

If we use Schwarz' lemma \cite{24}, we will get inequalities:
\begin{align}
\max_{1 \le i \le n}|[z^-(\Gamma_i^-)]^{(k-1)} |&\le q^-\max_{\Gamma^0}|[z(\Gamma^0)]^{(k-1)}|,\label{eq:53}\\
\max_{1 \le j \le m}|[z^+(\Gamma_j^+)]^{(k-1)} |&\le q^+\max_{\Gamma^0}|[z(\Gamma^0)]^{(k-1)}|, \label{eq:54}
\end{align}
where $q^+=const,\ 0<q^+<1,\ q^-=const,\ 0<q^-<1$. Note, that these
constants depend only on geometric properties of domains $D^-$ and
$D^+$.

If we use inequalities (\ref{eq:53}), (\ref{eq:54}), then we have
\begin{equation*}
\max_{\Gamma^0}|[z(\Gamma^0)]^{(k)}|\le q^+ \sum_{j=1}^m \beta_j^+ \cdot \max_{\Gamma^0}|[z(\Gamma^0)]^{(k-1)}|+ q^-\sum_{i=1}^n \beta_i^- \cdot \max_{\Gamma^0}|[z(\Gamma^0)]^{(k-1)}|, 
\end{equation*}
or
\begin{equation}\label{eq:55}
\max_{\Gamma^0}|[z(\Gamma^0)]^{(k)}|\le Q\max_{\Gamma^0}|[z(\Gamma^0)]^{(k-1)}|,
\end{equation}
where $Q=q^+ \sum_{j=1}^m \beta_j^+ + q^-\sum_{i=1}^n \beta_i^- .$

Taking into account the conditions $\beta_i^-=const> 0$,
$\beta_j^+=const> 0$, $0<\sum_{j=1}^m \beta_j^+ + \sum_{i=1}^n
\beta_i^- \le 1$, we obtain $0<Q<1$. This implies that
\begin{equation*}
  \lim_{k\to\infty}[z(\Gamma^0)]^{(k)}=0.
\end{equation*}

If the solution of the problem (\ref{eq:31})-(\ref{eq:37}) exists,
then by the maximum principle we obtain
\begin{align*}
  \max_{D^-}\left|[u^-(x_1,x_2 )]^{(k)}-u^-(x_1,x_2 )\right|&=O(Q^k),\\
  \max_{D^+}\left|[u^+(x_1,x_2 )]^{(k)}-u^+(x_1,x_2 )\right|&=O(Q^k),
\end{align*}
and, accordingly,
\begin{equation*}
  \max_{D}\left|[u(x_1,x_2 )]^{(k)}-u(x_1,x_2 )\right|=O(Q^k).
\end{equation*}
Thus, the iterative process (\ref{eq:43})-(\ref{eq:47}) converges to
this solution of the problem (\ref{eq:31})-(\ref{eq:37}) at the rate
of an infinitely decreasing geometric progression with ratio $Q$.

\end{proof}

\begin{Remark}
  By using the described iterative algorithm
  (\ref{eq:43})-(\ref{eq:47}) the solution of a non--classical contact
  problem (\ref{eq:31})-(\ref{eq:37}) is reduced to the solution of a
  sequence of classical boundary problems, which can be solved by any
  well-studied method.
\end{Remark}

\subsection{Existence of a Solution of Problem (\ref{eq:31})-(\ref{eq:37})}\label{sec:3.5}

Let us now prove the existence of a regular solution of the problem
(\ref{eq:31})-(\ref{eq:37}) in case of $f^-(x_1,x_2 )\equiv 0$ and
$f^+(x_1,x_2 )\equiv 0$. We introduce the notation
$\varepsilon^{(k)}(x_1,x_2)=u^{(k)}(x_1,x_2 )-u^{(k-1)}(x_1,x_2
)$. Then for the function $\varepsilon^{(k)}$ we obtain the following
problem
\begin{align*}
  L^-[\varepsilon^-(x_1,x_2)]^{(k)}=&0, && \text{ if }  (x_1,x_2 )\in D^-,
  \\
  L^+[\varepsilon^+(x_1,x_2)]^{(k)}=&0, && \text{ if } (x_1,x_2 )\in D^+,
  \\
  [\varepsilon^-(x_1,x_2 )]^{(k)}=&0, && \text{ if } (x_1,x_2 )\in \gamma_1^-\cup \gamma_2^-\cup \gamma_3,
  \\
  [\varepsilon^+ (x_1,x_2)]^{(k)}=&0, && \text{ if }  (x_1,x_2 )\in \gamma_1^+\cup \gamma_2^+\cup \gamma_4,
\end{align*}
\begin{equation*}
  [\varepsilon(\Gamma^0)]^{(k)}=[\varepsilon^- (\Gamma^0)]^{(k)}=[\varepsilon^+ (\Gamma^0)]^{(k)}
  =\sum_{i=1}^n \beta_i^- [\varepsilon ^-(\Gamma_i^-)]^{(k-1)}+\sum_{j=1}^m \beta_j^+ [\varepsilon^+(\Gamma_j^+)]^{(k-1)},     
\end{equation*}
where $k=1,2,\ldots,$ and
$[\varepsilon^-(\Gamma_i^-)]^{(0)}=0,$ $[\varepsilon^+
  (\Gamma_j^+)]^{(0)}=0,$
$i=\overline{(1,n)},\ j=\overline{(1,m)}.$

Then, analogously to (\ref{eq:55}), we obtain the estimation
\begin{equation*}
  \max_{\Gamma^0}|[\varepsilon(\Gamma^0)]^{(k)}|\le Q\max_{\Gamma^0}|[\varepsilon(\Gamma^0)]^{(k-1)}|,\qquad  0<Q<1,
\end{equation*}
or
\begin{equation*}
  \max_{\Gamma^0}|u^{(k)}(\Gamma^0)-u^{(k-1)}(\Gamma^0)|\le Q\max_{\Gamma^0}|u^{(k-1)}(\Gamma^0)-u^{(k-2)}(\Gamma^0)|,\qquad  0<Q<1,
\end{equation*}

This means that the sequence $\{u^{(k)}(x_1,x_2 )\}$ converges
uniformly on $\Gamma^0$. Then the functions $[u^+(x_1,x_2 )]^{(k)}$
and $[u^-(x_1,x_2 )]^{(k)}$ converge to the functions $u^+ (x_1,x_2 )$
and $u^- (x_1,x_2 )$, respectively, on the domains $D^-$ and $D^+$ on
the base of Harnack's first theorem \cite{25,26}.

From this we conclude that the limit function is the regular solution
of the problem (\ref{eq:31})-(\ref{eq:37}):

\begin{equation*}
  \lim_{k\rightarrow \infty}u^{(k)}(x_1,x_2)=u(x_1,x_2 ).
\end{equation*}

\section{Numerical Example}\label{sec:4}
\begin{figure}[htbp]
  \centering
  \includegraphics[width=0.5\linewidth]{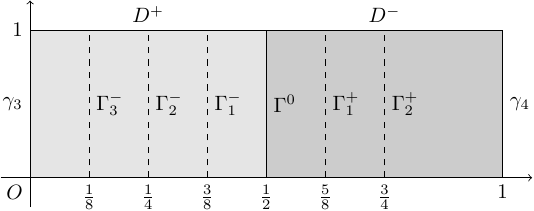}
  \caption{Domain for numerical example}
  \label{fig:3}
\end{figure}

Let us consider the area $D=\{(x_1,x_2 ) |0<x_1<1,\ 0<x_2<1\}$
(see the Figure \ref{fig:3}).

We consider the following test problem for the numerical solution:
find in $\bar{D}$ a continuous function (\ref{eq:31}) $u(x_1,x_2 )$,
which satisfies the following equations:

\begin{align*}
  \frac {\partial}{\partial x_1}[(1+x_1^2)\frac {\partial u^-}{\partial x_1}]+\frac {\partial}{\partial x_2}[(1+x_2^2)\frac {\partial u^-}{\partial x_2}]&=f^-(x_1,x_2), && \text{if  }(x_1,x_2)\in D^-,\\
\frac {\partial}{\partial x_1}[(1+2x_1^2)\frac {\partial u^+}{\partial x_1}]+\frac {\partial}{\partial x_2}[(1+2x_2^2)\frac {\partial u^+}{\partial x_2}]&=f^+(x_1,x_2), && \text{if  }(x_1,x_2)\in D^+,
\end{align*}
where
\begin{align*}
f^-(x_1,x_2)=&-\frac{1}{4}x_1 x_2(-16+\pi ^2(1+ x_2^2))\cos \frac{\pi x_2}{2}-\pi x_1(1+2 x_2^2)\sin \frac{\pi x_2}{2}\\
f^+(x_1,x_2)=&-4x_1 x_2 \cos \frac{\pi x_2}{2}\\
&+\frac{1}{4}(x_1-1)\left[x_2(-16+\pi ^2(1+2x_2^2))\cos \frac{\pi x_2}{2}+4\pi(1+4x_2^2)\sin \frac{\pi x_2}{2}\right],
\end{align*}

The function $u(x_1,x_2 )$ also satisfies the boundary conditions
\begin{align*}
  u^-(x_1,x_2 )&=0, && \text{if } (x_1,x_2 )\in \gamma_1^-\cup \gamma_2^-\cup \gamma_3,\\
  u^+(x_1,x_2)&=0, && \text{if } (x_1,x_2)\in \gamma_1^+\cup \gamma_2^+\cup \gamma_4,          
\end{align*}
the nonlocal contact conditions
\begin{multline*}
  u^-(\frac{1}{2},x_2)=u^+(\frac{1}{2},x_2)=  u(\Gamma^0) =
  \frac{1}{8} u^+\left(\frac{5}{8},x_2\right)+\frac{1}{8} u^+\left(\frac{3}{4},x_2\right)\\
  +\frac{1}{8} u^-\left(\frac{1}{8},x_2\right)+\frac{1}{8} u^-\left(\frac{1}{4},x_2\right)+\frac{1}{8} u^-\left(\frac{3}{8},x_2\right)
  +\frac {21}{64}x_2 \cos \frac{\pi x_2}{2},
\end{multline*}
and the coordination conditions are fulfilled.

The exact solution of this problem is
\begin{equation*} 
  u(x_1,x_2 )=
  \begin{cases}
    x_1 x_2\cos \frac{\pi x_2}{2}, & \text{ if } (x_1,x_2)\in D^-,\\
    \frac{1}{2} x_2\cos \frac{\pi x_2}{2}, & \text{ if } (x_1,x_2)\in \Gamma^0,\\
    (1-x_1)x_2\cos \frac{\pi x_2}{2}, & \text{ if } (x_1,x_2)\in D^+,
  \end{cases}  
\end{equation*}

Let us consider the following iterative process:
\begin{equation*}
  \frac {\partial}{\partial x_1}[(1+x_1^2)\frac {\partial (u^-)^{(k)}}{\partial x_1}]+\frac {\partial}{\partial x_2}[(1+x_2^2)\frac {\partial (u^-)^{(k)}}{\partial x_2}]=f^-(x_1,x_2),\ \mbox{if  }(x_1,x_2)\in D^-,
\end{equation*}
\begin{equation*}
  \frac {\partial}{\partial x_1}[(1+2x_1^2)\frac {\partial (u^+)^{(k)}}{\partial x_1}]+\frac {\partial}{\partial x_2}[(1+2x_2^2)\frac {\partial (u^+)^{(k)}}{\partial x_2}]=f^+(x_1,x_2),\ \mbox{if  }(x_1,x_2)\in D^+,
\end{equation*}
\begin{equation*}
  [u^-(x_1,x_2 )]^{(k)}=0, \  \mbox{if } (x_1,x_2 )\in \gamma_1^-\cup \gamma_2^-\cup \gamma_3,
\end{equation*}
\begin{equation*}
  [u^+(x_1,x_2)]^{(k)}=0, \  \mbox{if } (x_1,x_2)\in \gamma_1^+\cup \gamma_2^+\cup \gamma_4,          
\end{equation*}
\begin{multline*}
  [u(\Gamma^0)]^{(k)}
  = \left[u^-\left(\frac{1}{2},x_2\right)\right]^{(k)}
  =\left[u^+\left(\frac{1}{2},x_2\right)\right]^{(k)}\\
  =\frac{1}{8} \left[u^-\left(\frac{1}{8},x_2\right)\right]^{(k-1)}
  +\frac{1}{8} \left[u^-\left(\frac{1}{4},x_2\right)\right]^{(k-1)}
  +\frac{1}{8} \left[u^-\left(\frac{3}{8},x_2\right)\right]^{(k-1)}\\
  +\frac{1}{8} \left[u^+\left(\frac{5}{8},x_2\right)\right]^{(k-1)}
  +\frac{1}{8} \left[u^+\left(\frac{3}{4},x_2\right)\right]^{(k-1)}
  +\frac {21}{64}x_2 \cos \frac{\pi x_2}{2},
\end{multline*}
where $k=1,2,\ldots,$ and the initial value
$[u(\Gamma^0)]^{(1)}$ is equal to $0$.

\begin{figure}[htbp]
  \centering
  \includegraphics[height=.3\linewidth]{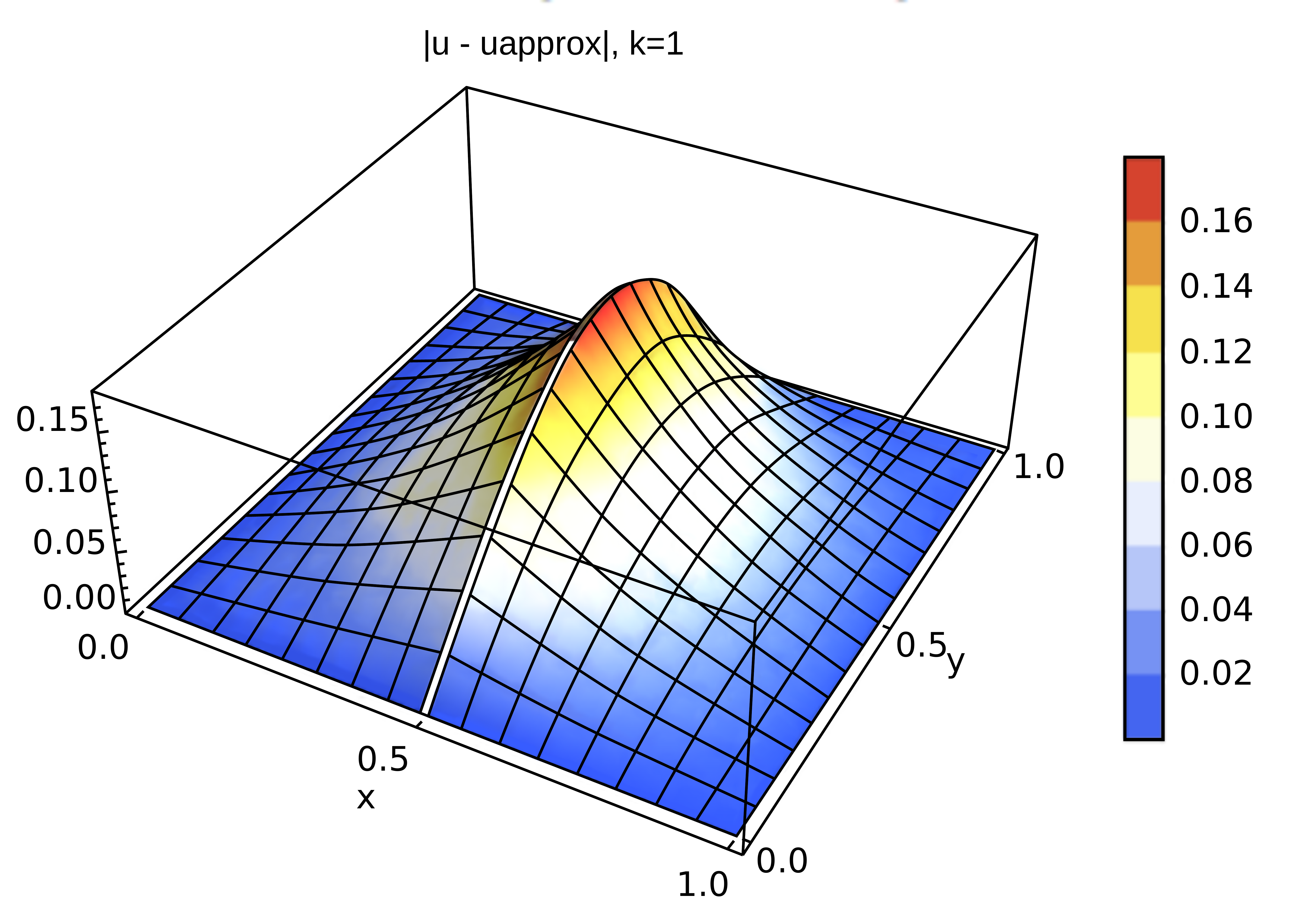}
  \includegraphics[height=.3\linewidth]{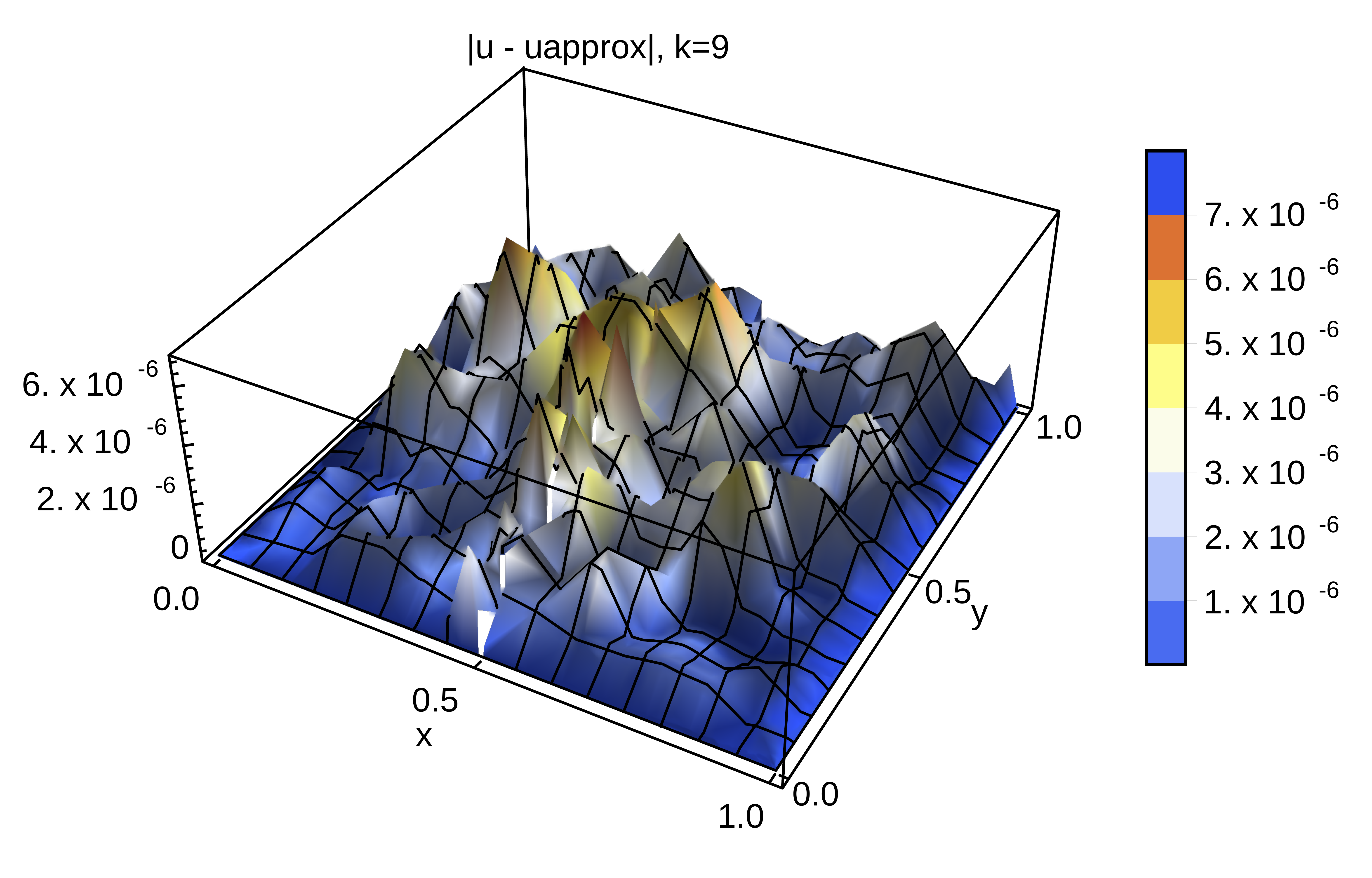}
  \caption{The magnitude of the difference between exact and
    approximate solutions for k = 1 and 9. For k=1, range is
    $0.02-0.16$; for k=9, range is $1.0\times 10^{-6} - 6.0\times
    10^{-6}.$}
  \label{fig:4}
\end{figure}

The numerical calculations were carried out using the program Wolfram
Mathematica (see Figure \ref{fig:4}). In computations the following
Wolfram functions were used with the corresponding arguments:
NDSolveValue (to assign the computed value to each component of the
vector-function -- the solution of the system of two equations with
partial derivatives), DirichletCondition (to specify boundary values
within the function NDSolveValue), Piecewise (to represent the
solution on the whole area for further visualization), NMaxValue (to
calculate uniform norm), Do (to organize the outer loop for
iterations), along with other supplementary functions to calculate
norms and get respective graphs.

From the Theorem \ref{thm:3}, we have $Q = q^+\sum_{j=1}^m \beta_j^+ +
q^- \sum_{i=1}^n \beta_i^-$, where $0 < q^+, q^- < 1$.  The absolute
error decreases approximately as $O(Q^k)$. Considering the nonlocal
contact condition of the numerical example, $Q < 5/8$.

\begin{figure}[htbp]
  \centering
  \includegraphics[width=0.9\linewidth]{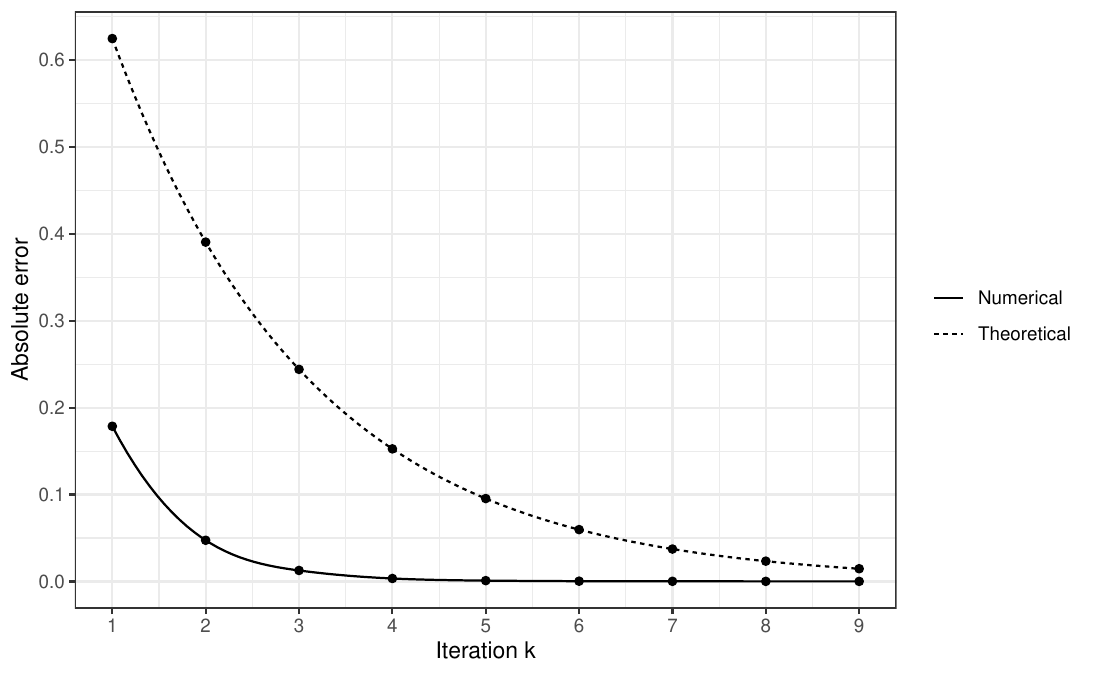}
  \caption{Absolute error (numerical) and $Q^k$ versus iteration $k$ }
  \label{fig:5}
\end{figure}

Figure \ref{fig:5} compares the absolute error (in C-norm) with
theoretical value $\left(\frac{5}{8}\right)^k$.

\begin{table}[htbp]
  \begin{center}
    \includegraphics{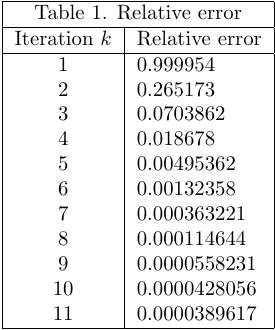}
  \end{center}
  \caption{Relative error}\label{tbl:1}
\end{table}

\begin{figure}[H]
   \begin{minipage}{0.45\textwidth}
     \centering
     \includegraphics[width=\linewidth]{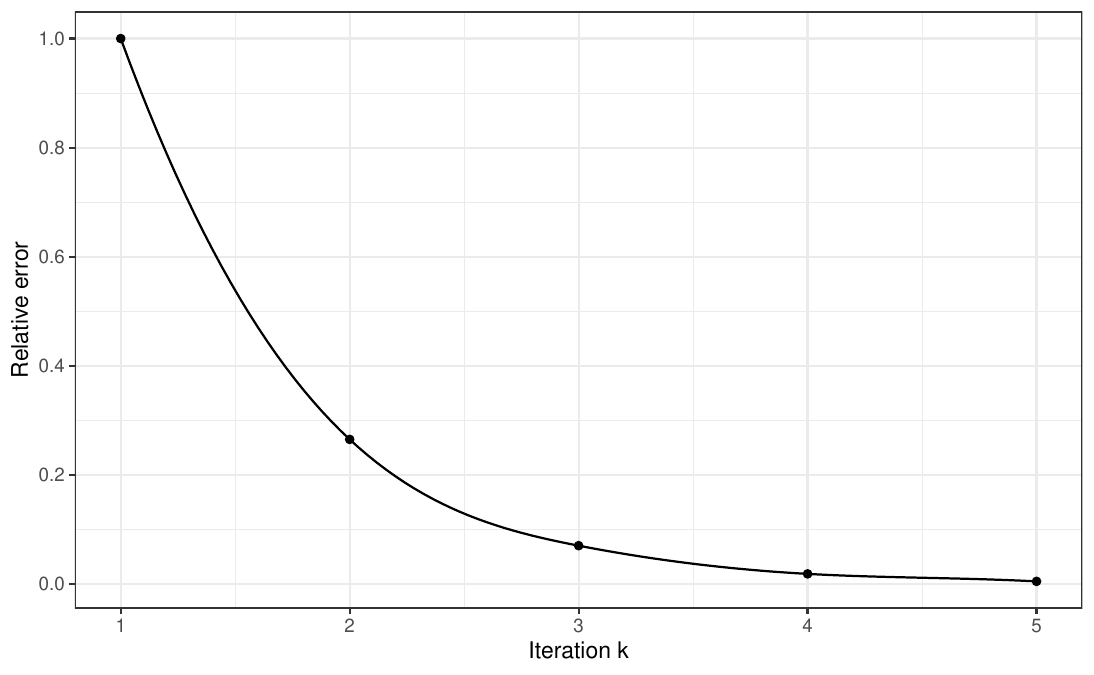}
     \caption{Relative error versus iteration k, k=1..5}\label{fig:6}
   \end{minipage}\hfill
   \begin{minipage}{0.45\textwidth}
     \centering
     \includegraphics[width=\linewidth]{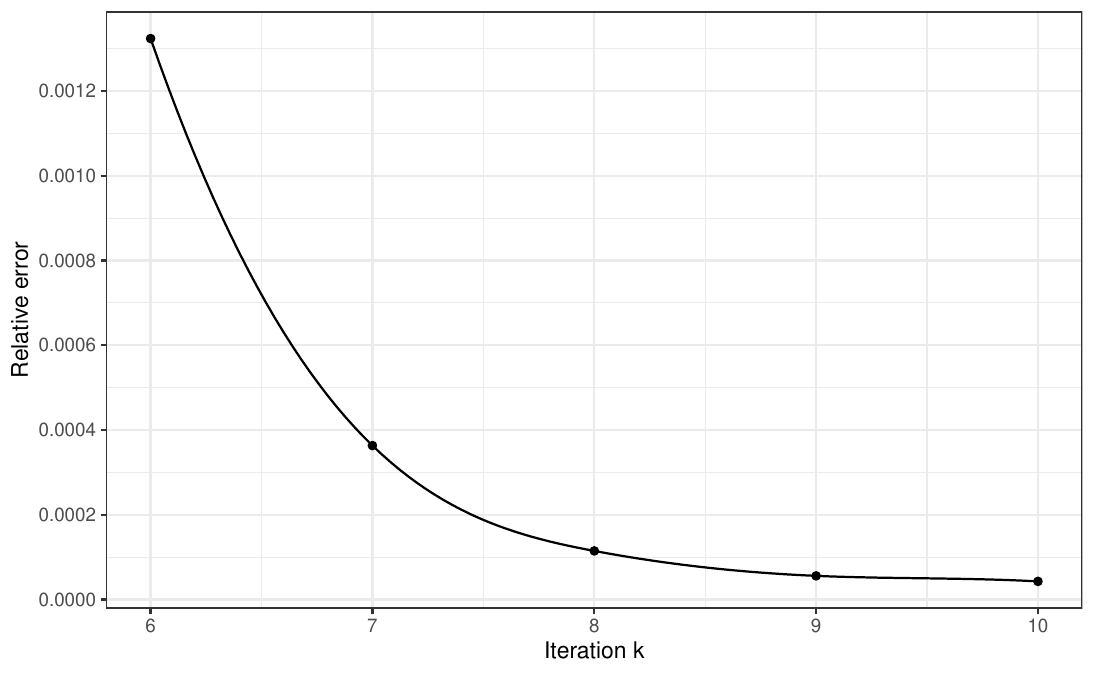}
     \caption{Relative error versus iteration k, k=6..10}\label{fig:7}
   \end{minipage}
\end{figure}

Figures \ref{fig:6}-\ref{fig:7} show the behavior of relative error
$\frac{\|u_{exact}-u_{appr} \|_C}{\|u_{exact} \|_C}$ (uniform C-norm
is taken on the open area $D$) versus iteration $k$, for
$k=1,\ldots,5$ and $k=6,\ldots,10$, respectively (see also results in
Table \ref{tbl:1}).

The results of numerical calculations fully agree with the theoretical
conclusions and show the efficiency of the proposed iterative
procedure.

\section{Conclusion}\label{sec:5}

The theory of contact problems is widely used in many fields of
mechanics (including construction mechanics), in mechanical
engineering, etc. In these problems, various contact conditions are
considered along the contact line (see, for example, \cite{27,28}).

In the present article a new type of nonclassical boundary-value
problem with nonlocal contact conditions along the contact line is
considered for elliptic equation with variable coefficients and mixed
derivatives. Thus, using the results of the present article, one can
expand the class of contact problems.

The main results of the proposed article can be formulated as follows:

The existence and uniqueness of the solution of a nonlocal contact
problem for the elliptic equation with variable coefficients and mixed
derivatives is proved. For this aim the convergent iterative method
(\ref{eq:43})-(\ref{eq:47}) is constructed, which also is used to find
the numerical solution. The method converges to the solution of the
problem (\ref{eq:31})-(\ref{eq:37}) at the rate of infinitely
decreasing geometric progression. By using this iterative algorithm
the solution of a non-classical contact problem is reduced to the
solution of a sequence of classical boundary problems, which can be
solved applying well-studied methods. The results of numerical
calsulations agree with theoretical results.

The analytical solution in a form of series is received for the same
problem, but with constant coefficients to avoid huge
formulae. Moreover, the applied technique can be successfully used for
more general problems, but in this case the use of spectral theory of
linear operators will be needed.

In contrast to boundary nonlocal problems, convergence is achieved
under the more general conditions: $0 < \gamma^-+\gamma^+ \le 1$ (in
the case of the Fourier method) and $0 < \sum_{j=1}^m \beta_j^+ +
\sum_{j=1}^n \beta_i^- \le 1$ (in the case of general equation with
variable coefficients).

The technique used in the present article can also be applied for the
problems with parabolic type equations.

%
%
%
%
%
%
%
%

\section{Declarations}

\paragraph{Availability of data and materials}\ 

\noindent Data sharing is not applicable to this article as no
datasets were generated or analysed during the current study.

\paragraph{Competing interests}\ 

\noindent The authors declare that they have no competing interests.

\paragraph{Funding}\ 

\noindent This research received funding from the University of Malaga
for its open access publication.

\paragraph{Authors' contribution}\ 

\noindent All authors contributed to all sections. All authors read
and approved the final manuscript.

\paragraph{Acknowledgements}\ 

\noindent Not applicable.


\providecommand{\bysame}{\leavevmode\hbox to3em{\hrulefill}\thinspace}
\providecommand{\MR}{\relax\ifhmode\unskip\space\fi MR }
\providecommand{\MRhref}[2]{%
  \href{http://www.ams.org/mathscinet-getitem?mr=#1}{#2}
}
\providecommand{\href}[2]{#2}

\end{document}